\documentclass[12pt]{article}
\usepackage{amssymb,amsmath,amsfonts,t1enc}
\newtheorem{tw}{\bf\indent Theorem.}
\newcommand{\card}{{\rm card\,}}

\begin{document}
\baselineskip=12truept
\centerline{\normalsize \bf Subobjects of the Successive
Power Objects in the Topos $G-Set$}
\vskip 1.0truecm
\centerline{\normalsize \bf Apoloniusz Tyszka}
\vskip 1.0truecm
\par
\noindent
{\bf Abstract.} Let $G$ be a group and let
$M$ be an object of the topos $G-Set.$ We prove that an
object $X$ of the category $G-Set$ is isomorphic to some
subobject of one of the objects $P(M),$ $P(P(M))$,
$P(P(P(M))),\dots$ if and only if card $X<\sup\{\card
P(M),\card P(P(M)),\card P(P(P(M))),\dots\}$ and $\{g\in
G:\forall m\in M\ gm=m\}\subseteq\{g\in G:\forall x\in X\
gx=x\}.$
\par
\leftskip=0.0truein
\rightskip=0.0truein
\baselineskip=14truept
\vskip 1.0truecm
In this paper we present a simplified proof of a theorem
which was proved in \cite{Tyszka1993}; earlier weaker results of this
type can be found in \cite{Szocinski} and \cite{Tyszka1992}.
Let $G$ be a group. The
category $G-Set$ has as objects, sets equipped with a
$G$-action and as morphisms functions preserving this
action. This is a topos (see \cite{Goldblatt} for the introduction to
the topoi theory), i.e. $G-Set$ is an elementary topos in
the terminology of \cite{Johnstone}. Let $P:G-Set \rightarrow G-Set$
denote the power functor.
\begin{tw}
Let $M$ be an object of the topos $G-Set$. Then an object
$X$ of the topos $G-Set$ is isomorphic to some subobject of
one of the objects $M_n:=\underbrace{P(\dots(P}_{n\ times}(M)))$
$(n \geq 1)$ if and only if
\begin{equation}
\card X<\chi(M):=\sup\{\card M, \card P(M), \card P(P(M)),\dots\},
\end{equation}
and
\begin{equation}
G_M:=\{g\in G:\forall m\in M\ gm=m\}\subseteq
G_X:=\{g\in G:\forall x\in X\ gx=x\}.
\end{equation}
\end{tw}
\vskip 0.2truecm
\par
\noindent
{\it Proof.} The necessity is obvious, we show the
sufficiency. Each transitive $G$-set $X$ which satisfies (2)
is isomorphic to a $G$-set of the form $(\{gH:g \in
G\},G,L)$ where $H$ is a subgroup of $G$ and $L$ denotes
left translation (see \cite{Bourbaki} p.106), obviously $G_M \subseteq
G_X\subseteq H.$ Take any relation which well orders $M$,
let $M=\{m_\gamma\}_{\gamma<\alpha}.$ Then the transitive
$G$-set determined by
$z:=\{\{m_\gamma:\gamma<\beta\}:\beta\leq\alpha\}\in M_2$
is isomorphic to $(\{gG_M:g\in G\}, G, L)$, hence
(Z.~Moszner in \cite{Tyszka1993}) $(\{gHz:g\in G\},G,L)$ i.e.
$(\{\{ghz: h \in H\}:g \in G\},G, L)$ is a subobject of $M_3$
which is isomorphic to $(\{gH:g \in G\},G,L)$.
\par
An arbitrary $G$-set $X$ which satisfies (1) and (2) is a
disjoint union of less than $\chi(M)$ transitive $G$-sets
$Y_\delta$ such that $G_M\subseteq G_{Y_\delta}:=\{g \in
G:\forall y \in Y_\delta \ gy=y\}.$ The proof will be
completed by showing that $M_{n+4}$ includes card $M_n$
pairwise disjoint objects, each of them isomorphic to
$M_3.$ Let $Bij(M)$ denote the group of bijections of $M$,
by ${\cal P}:Set \rightarrow Set$ we mean the power set
functor. The relation:
$$
xR_{n+1}y:\Longleftrightarrow \exists f \in Bij(M)\ y=\underbrace{{\cal
P}(\dots({\cal P}}_{n+1\ times}(f)))(x),\ \ \mbox{where}\
x,y\in M_{n+1}
$$
is an equivalence relation on $M_{n+1}.$ Let
$[x]_{R_{n+1}}$ denote the equivalence class of $x \in
M_{n+1},$ obviously card $[x]_{R_{n+1}}\leq \card\
Bij(M).$ Using the inequality $\card\ M_{n+1}>\card
M_n \cdot \card \ Bij(M)$ (valid both for finite and
infinite $M$) we obtain that
\linebreak
$\card\ (M_{n+1}/
R_{n+1})>\card\ M_n.$ If $w\in M_{n+1}$ then $\forall
g\in G\ g[w]_{R_{n+1}}=[w]_{R_{n+1}},$ hence the following
$G$-set:
$$
X_w:=\left\{\langle\underbrace{\{\dots\{x\}\dots\}}_{n-1\
times\{\ \}}, [w]_{R_{n+1}}\rangle:x\in M_3\right\}
$$
is a subobject of $M_{n+4}$ which is isomorphic to $M_3$.
Moreover, for every $s, t\in M_{n+1}$ if $[s]_{R_{n+1}} \neq
[t]_{R_{n+1}}$ then $X_s\cap X_t=\emptyset.$ This
observation completes the proof.

Apoloniusz Tyszka\\
Technical Faculty\\
Hugo Ko\l{}\l{}\k{a}taj University\\
Balicka 116B, 30-149 Krak\'ow, Poland\\
E-mail address: {\it rttyszka@cyf-kr.edu.pl}

\begin{thebibliography}{6}
\bibitem{Bourbaki}
N. Bourbaki, {\it \'El\'ements de Math\'ematique,
Livre II, Alg\`ebre, Chapitre I, Struc\-tu\-res alg\'ebriques,}
Hermann, Paris 1958.
\bibitem{Goldblatt}
R. Goldblatt, {\it Topoi, the Categorial Analysis of
Logic,} North-Holland, Amsterdam 1979.
\bibitem{Johnstone}
P. T. Johnstone, {\it Topos Theory,} Academic Press,
New York 1977.
\bibitem{Szocinski}
B. Szoci\'nski, {\it Basic concepts of Klein
geometries,} Zesz. Nauk. Politech. \'Sl. 1055, Mat.-Fiz. 62, 1990. 
\bibitem{Tyszka1992}
A. Tyszka, {\it On one combinatorial lemma and its
geometric consequences},\\ Demonstr. Math. 25, 1992, pp. 579--582.
\bibitem{Tyszka1993}
A. Tyszka, {\it On the notion a geometric object in a
Klein space,} Rocznik Naukowo-Dydaktyczny Wy\.zszej Szko\l{}y
Pedagogicznej w Krakowie 159, Prace Matematyczne XIII
(Annales de l'Ecole Normale Sup\'erieure \`a Cracovie 159,
Travaux Mathematiques XIII), Krak\'ow 1993, pp. 287--299.
\end{thebibliography}
\end{document}